\documentclass[12pt,reqno]{article}
\def\hybrid{\topmargin      0pt
\oddsidemargin 0pt \headheight 0pt \headsep 0pt \textwidth 165true
mm
\textheight 235true mm
\marginparwidth 0.0in \parskip  0pt plus 1pt \jot = 1.5ex}

\usepackage{amssymb}
\usepackage{amsmath}
\usepackage{amsthm}
\hybrid
\usepackage{amssymb}
\usepackage{amsmath}
\usepackage{amsthm}
\newcommand{\be}[1]{\begin{eqnarray#1}}
\newcommand{\ee}[1]{\end{eqnarray#1}}
\newtheorem{thm}{Theorem}[section]
\newtheorem{propn}[thm]{Proposition}
\newtheorem{lemma}[thm]{Lemma}
\newtheorem{corollary}[thm]{Corollary}

\theoremstyle{definition}
\newtheorem{remark}[thm]{Remark}

\newcommand{\tp}{\otimes}
\newcommand{\ot}{\otimes}

\newcommand{\braid}[2]{{#1}$\lower4pt\hbox{$\tp\atop\raise4pt
            \hbox{$\scriptscriptstyle\Ru $}$}${#2}}
\newcommand{\twist}[2]{{#1}${\,\scriptscriptstyle
\Ru}\atop\raise9pt\hbox{$\scriptstyle\tp$} ${#2}}
\newcommand{\twistF}[2]{{#1}${\,\scriptscriptstyle
\F}\atop\raise9pt\hbox{$\scriptstyle\tp$} ${#2}}

\newcommand{\J}{\mathcal{J}}
\newcommand{\A}{\mathcal{A}}
\renewcommand{\L}{\mathcal{L}}
\newcommand{\I}{\mathfrak{I}}
\newcommand{\Ec}{\mathcal{E}}

\renewcommand{\O}{O}

\renewcommand{\S}{\mathcal{S}}

\newcommand{\Ru}{\mathcal{R}}
\newcommand{\V}{\mathcal{V}}

\newcommand{\U}{\mathcal{U}}
\newcommand{\F}{\mathcal{F}}
\newcommand{\C}{\mathbb{C}}
\newcommand{\Z}{\mathbb{Z}}

\newcommand{\N}{\mathbb{N}}

\newcommand{\Ind}{\mathrm{I}}

\newcommand{\ff}{\varphi}

\newcommand{\QQ}{\mathcal{Q}}

\newcommand{\LL}{\mathcal{L}}
\newcommand{\KK}{\mathcal{K}}

\newcommand{\g}{\mathfrak{g}}
\newcommand{\p}{\mathfrak{p}}
\newcommand{\n}{\mathfrak{n}}

\renewcommand{\l}{\mathfrak{l}}
\renewcommand{\c}{\mathfrak{c}}

\newcommand{\la}{\lambda}

\newcommand{\nn}{\nonumber}
\newcommand{\al}{\alpha}

\newcommand{\id}{\mbox{id}}
\newcommand{\End}{\mathrm{End}}

\newcommand{\Char}{\mathrm{ch }}

\newcommand{\Tr}{\mathrm{Tr}}
\newcommand{\T}{\mathrm{T}}
\newcommand{\si}{\sigma}

\newcommand{\h}{\mathfrak h}

\begin{document}
\title{Explicit equivariant quantization on coadjoint orbits of $GL(n,\C)$}
\author{J. Donin$^{\dag,\ddag,}$\footnote{This research is partially supported
by the Israel Academy of Sciences grant no. 8007/99-01.}
\hspace{3pt} and A. Mudrov$^\dag$}
\date{}
\maketitle
\begin{center}
{\dag Department of Mathematics, Bar Ilan University, 52900 Ramat Gan,
Israel.\\
\ddag Max-Plank-Institut. f\"ur Mathematik, Vivatsgasse 7,
D-53111 Bonn. }
\end{center}
\begin{abstract}
We present an explicit $\U_h\bigl(gl(n,\C)\bigr)$-equivariant
quantization on coadjoint orbits of $GL(n,\C)$. It forms a
two-parameter family quantizing the Poisson pair of the reflection
equation and Kirillov-Kostant-Souriau brackets.
\end{abstract}

\section{Introduction}

Let $G$ be a reductive complex connected Lie group with Lie
algebra $\g$. Let $\U(\g)$ be the universal enveloping algebra of
$\g$ and $\U_h(\g)$ the corresponding quantum group, \cite{Dr}. We
consider the problem of $\U_h(\g)$-equivariant quantization on
$\g^*$ and on coadjoint orbits of $G$.

By a $\U_h(\g)$-equivariant quantization on a $G$-manifold $M$ we
mean a two-parameter family  $\A_{t,h}$ of  $\U_h(\g)$-module
algebras such that $\A_{0,0}=\A(M)$, the function algebra on $M$,
and the $\U_h(\g)$-action on $\A_{t,h}$ is an extension of the
$\U(\g)$-action on $\A(M)$. Briefly, we refer to such a
quantization as $(t,h)$-quantization while to its $(0,h)$- and
$(t,0)$-subfamilies as $h$- and $t$-quantizations. This definition
implies that the family $\A_{t,0}$ is $G$-equivariant.

At the infinitesimal level, a $(t,h)$-quantization gives rise to a
pair of compatible Poisson brackets.  The one along the parameter
$h$ is represented as a difference of bivector fields, $f-r_M$,
\cite{D2}, where $f$ is a $G$-invariant one and $r_M$ is induced
on $M$  via the group action by the r-matrix $r\in \wedge^2 \g$ corresponding to
$\U_h(\g)$. The Poisson bracket along the
parameter $t$ is $G$-invariant.

In the paper, as $G$-manifolds we consider $\g^*$ and coadjoint
orbits, as a $G$-invariant Poisson bracket along $t$ we take the
Poisson-Lie bracket on $\g^*$ and its restriction to an orbit; the
latter is known as the Kirillov-Kostant-Souriau (KKS) bracket. As
was shown in \cite{D1}, a $(t,h)$-quantization on $\g^*$ exists
only for $\g$ with simple components being $sl(n,\C)$. So we
assume $\g=gl(n,\C)$  in the paper.

Firstly, we show that the quantization on $\g^*$ is given by the
so-called extended reflection equation (ERE) algebra, $\L_{t,h}$.
At $(t,h)=(0,0)$, it is the polynomial algebra on $\g^*$ (or
symmetric algebra of $\g$). The ordinary or quadratic reflection
equation (RE) algebra $\L_{h}$ was studied in \cite{KSkl,KS}. The
ERE algebra naturally appears in non-commutative differential
calculus related to quantum groups, see e.g. \cite{IP}.

Secondly, using the generalized Verma  module over $\U(\g)$, we
show that the ERE algebra can be restricted to any semisimple
orbit in $\g^*$. By the restriction we mean the following. Let
$\{f\}\subset \A(\g^*)$ be a set of functions generating the ideal
of the orbit $O$. Then, there exist their extensions
$\{f_{t,h}\}\subset\L_{t,h}$ generating an ideal in $\L_{t,h}$;
the quotient of $\L_{t,h}$ by that ideal gives a flat deformation
over $\C[[t,h]]$.

Thirdly, we give explicit formulas for the generators
$\{f_{t,h}\}$.

Let us describe our approach in more detail. In the classical
case, an orbit of rank $k-1$ is specified by a matrix polynomial
equation, \be{}\label{mpe1} (X-\mu_1)\ldots(X-\mu_k)=0, \ee{} and
the conditions on traces
$$
\Tr(X^m)=\sum_{i=1}^k n_i \mu^m_i, \quad m=1,\ldots,k-1,
$$
where $n_i$ are non-negative integers such that $\sum_{i=1}^k n_i
=n$. This orbit consists of matrices with eigenvalues $\mu_i$ of
multiplicities $n_i$, $i=1,\ldots,k$. The ideal of the orbit is
generated by the $n\times n$ entries of the matrix polynomial and
$k-1$ functions involving traces. Matrix polynomial equation
(\ref{mpe1}) defines a $G$-invariant subvariety $M^k_\mu$, which
is a finite collection of rank $<k$ orbits at generic $\mu$,
namely, when $\mu_i\not = \mu_j$ if $i\not = j$. The condition on
traces defines a character of the subalgebra  $\S$ of invariants
restricted to $M^k_\mu$. This restriction is a finite dimensional
semisimple algebra and it has a finite spectrum whose points
correspond to orbits.

We extend this picture to the quantum case. It turns out that the
quotient of the algebra $\L_{t,h}$ by the relations
$$
(L-\mu_1)\ldots(L-\mu_k)=0,
$$
where $L$ is the matrix whose entries are generators of
$\L_{t,h}$, gives a $(t,h)$-quantization, $\A_{t,h}(M^k_\mu)$, on
$M^k_\mu$. Moreover, the quotient of $\L_{t,h}$ by the ideal
generated by entries of {\em any} polynomial in $L$ is a flat
algebra over $\C[[t,h]]$. The first part of the paper is devoted
to a proof of this statement.

With the matrix polynomial condition imposed, the equations on
traces have to be deformed in a consistent way. The quantum
version of these equations involves the quantum trace $\Tr_q$,
where $q=e^h$, that is the quantum analog of the ordinary trace,
\cite{FRT}. It is defined in such a way that the elements $\Tr_q(L^m)$,
$m\in \N$, are $\U_h(\g)$-invariant. The first $n$ traces
$\Tr_q(L^m)$, $m=1,\ldots,n$, generate the center of the algebra
$\L_{t,h}$. The problem of restricting the quantization
$\A_{t,h}(M^k_\mu)$ to particular orbits in $M^k_\mu$ reduces to
the problem of computing the values of quantum traces on them. The
second part of the paper is devoted to solution of this problem.

Let us illustrate our results on the simplest example of symmetric
(rank one) orbits. The quantized orbit of matrices with
eigenvalues $\mu_1,\mu_2$ of multiplicities $n_1, n_2$ is the
quotient of the ERE algebra $\L_{t,h}$ by the relations \be{}
(L-\mu_1)(L-\mu_2)&=&0,
\nn\\
\Tr_q(L)&=&\hat n_1\mu_1+\hat n_2\mu_2+t \hat n_1\hat n_2. \nn
\ee{} In the paper, we define $\hat m$ as
$\frac{1-q^{-2m}}{1-q^{-2}}$ for $m\in \Z$ and normalize the
quantum trace so that $\Tr_q(1)=\hat n$. Note that these
particular formulas for symmetric orbits were derived in our paper
\cite{DM2} by different methods.

As a limit case of our quantization, we come, putting $h=0$, to an
explicit quantization of the KKS bracket. The problem of
quantizing  this Poisson structure was put forward in
\cite{BFFLS}. It is interesting that the $G$-invariant
quantization of the KKS bracket is obtained with the substantial
use of quantum groups.

Among other works relevant to our study, we would like to mention
\cite{Ast,Kar}, where the Karabegov quantization with separation
of variables is applied to building the quantum moment
map. That construction establishes relations between
$G$-equivariant quantizations and generalized Verma modules. A
construction of the $G$-equivariant quantization on semisimple
orbits via the generalized Verma modules is presented in
\cite{DGS}.

None of the mentioned methods gives explicit formulas for
quantized orbits. In \cite{DGK}, explicit formulas was given for
quantizing $\C \mathrm{P}^n$ type orbits. In our papers \cite{DM1,DM2}, we
developed a method of quantization on orbits of small ranks and
built an explicit quantization on symmetric and bisymmetric orbits
of $gl(n,\C)$. In the present paper, we generalize those formulas
to all semisimple orbits.

Let us remark that our quantized orbits are relevant to the so-called fuzzy
spaces considered in the physics literature (see, for example, \cite{DolJ,PawSt}
and references therein).

The paper is organized as follows. Section \ref{sKKS} deals with
the $t$-quantization of the KKS bracket. There, we prove that the
natural $G$-equivariant quantization on $\g^*$ can be restricted
to semisimple orbits. In Section \ref{sQGC} we generalize that
approach to the $(t,h)$-quantization. Therein, we prove that the
quotient of $\L_{t,h}$ by any matrix polynomial equation is a flat
deformation over $\C[[t,h]]$. In Section \ref{sInv},
we study the subalgebra of invariants in $\L_{t,h}$ restricted to
$M^k_\mu$. We prove that it is a finite dimensional semisimple
commutative algebra and compute its characters. The explicit
formulas for the equivariant quantization on orbits are given in
Section \ref{sQO}.

\vspace{0.2in}

\noindent {\bf \large Acknowledgment.} We are grateful to J.
Bernstein and S. Shnider for  helpful discussions and valuable
comments.
\section{$G$-equivariant quantization on orbits}
\label{sKKS}
\subsection{Quantization on $\g^*$}

Let $\U(\g)$ be the universal enveloping algebra of a complex Lie
algebra $\g$ with the Lie bracket
$[\hspace{3.pt}.\hspace{2.pt},.\hspace{2.5pt}]$. Let $G$ be a
connected Lie group corresponding to $\g$. We consider $\g[t]$ as
a Lie algebra over $\C[t]$ with respect to the bracket
$[x,y]_t=t[x,y]$. Let $\L_t$ denote the universal enveloping
algebra of $\g[t]$. By definition, it is a quotient of the tensor
algebra $\T(\g)[t]$ by the ideal generated by relations $x\ot
y-y\ot x-t[x,y]$, $x,y\in \g$. The algebra $\L_t$ is a
$G$-equivariant quantization of the symmetric algebra of $\g$
(considered as the polynomial algebra on $\g^*$) with the
Poisson-Lie bracket induced by
$[\hspace{3.pt}.\hspace{2.pt},.\hspace{2.5pt}]$. The assignment
$x\mapsto tx$, $x\in\g$, defines a $\C[t]$-algebra morphism,
\be{}\label{mult} \phi_t\colon\L_t\to U(\g)[t], \ee{} which is
obviously an embedding of free modules over $\C[t]$.

\subsection{Levi and parabolic subalgebras}

Let $\g$ be a complex reductive Lie algebra with the Cartan
subalgebra $\h$ and $\g=\n^-\oplus\h\oplus\n^+$ its polarization
with respect to $\h$.

We fix a Levi subalgebra $\l$, which is, by definition, the
centralizer of an element in $\h$. The algebra $\l$ is reductive,
so it is decomposed into the direct sum of its center and the
semisimple part, $\l=\c\oplus [\l,\l]$. Also, there exists a
decomposition \be{} \label{levi} \g=\n_\l^-\oplus\l\oplus\n_\l^+,
\ee{} where $\n_\l^\pm$ are subalgebras in $\n^\pm$.

It is clear that $\c\subset\h$ and $\h=\c\oplus\h_\l$, where
$\h_\l=\h\cap[\l,\l]$. So, we have a projection $\h\to\c$. On the
other hand, the Cartan decomposition defines a projection \be{}
\label{projection} \pi\colon\g\to\h. \ee{} Composition of these
maps gives the natural projection $\g\to\c$ which defines an
embedding $\c^*\to\g^*$.

Taking an element $\mu\in\c^*$, we consider the coadjoint orbit of
$G$ in $\g^*$ passing through $\mu$. For generic $\mu$, this orbit
is semisimple with $\l$ being the Lie algebra of the stabilizer at
the point $\mu$. Denote by $\O^k$, $k=\dim \c$, the closure in
$\g^*\times\c^*$ of this family of semisimple orbits and by
$\O^k_\mu$ the fiber over $\mu$. It is known that the family
$\O^k$ is flat over $\c^*$. For generic $\mu\in\c^*$ the variety
$\O^k_\mu$ is a semisimple orbit in $\g^*$.

Let $\p$ denote the parabolic subalgebra $\l\oplus\n_\l^+$.
Decomposition (\ref{levi}) turns into the decomposition \be{}
\label{parabolic} \g=\n_\l^-\oplus\p. \ee{} It is clear that
characters (one-dimensional representations) of $\p$ (as well as
of $\l$) are generated by linear forms from $\c^*$.

\subsection{Generalized Verma modules}

For any $\mu\in\c^*$, projection (\ref{projection}) defines a
representation $\pi_\mu$ of $\p$ on $\C$:
$\pi_\mu(x)=\mu\bigl(\pi(x)\bigr)$, $x\in\p$. It extends to a
representation of $\U(\p)$, which we still denote by $\pi_\mu$.

The generalized Verma $\U(\g)$-module corresponding to
$\mu\in\c^*$ is the left module $\V_\mu=\U(\g)\ot_{\U(\p)}\C$. In
the tensor product, $\C$ is a $\U(\p)$-module with respect to
$\pi_\mu$, while $\U(\g)$ is considered as a right
$\U(\p)$-module. The element $v_0=1\ot 1\in \V_\mu$ is the highest
weight vector. In particular, $xv_0=\pi_\mu(x)v_0$ for $x\in
\U(\p)$. Moreover, the map $\U(\n_\l^-)\to \V_\mu$, $x\mapsto
xv_0$, is an isomorphism of vector spaces. So, all $\V_\mu$ are
canonically isomorphic to $\V=\U(\n_\l^-)$ as vector spaces. The
representation $\V_\mu$ is given by a homomorphism of algebras
\be{} \label{ff} \ff_\mu\colon\U(\g)\to \End(\V). \ee{} We treat
$\V_\mu$ as belonging to the family of Verma modules parameterized
by $\c^*$. Namely, we consider trivial bundles over $\c^*$ with
the fibers $\U(\g)$ and $\End(\V)$; then $\ff_\mu$ is a map of
their polynomial sections, which we denote by $\U(\g)[\mu]$ and
$\End(\V)[\mu]$. One can think of $\mu$ as a collection of formal
coordinates on $\c^*$, $\mu=(\mu_1,...,\mu_k)$. Then,
$\U(\g)[\mu]$ and $\End(\V)[\mu]$ are polynomials in $\mu$ with
values in $\U(\g)$ and $\End(\V)$, respectively.

It is easy to check that the map
$\ff_{\mu/t}\circ\phi_t\colon\L_t\to \End(\V)[[t]]$, where
$\phi_t$ is given by (\ref{mult}), defines a map of polynomial
sections over $\c^*$, \be{}\label{Polsec}
\Phi_{t,\mu}\colon\L_t[\mu]\to \End(\V)[[t]][\mu]. \ee{}

\begin{thm}[\cite{Ast,DGS}]
\label{thm1} The image $\A_{t,\mu}$ of $\L_{t,\mu}$ in
$\End(\V)[[t]][\mu]$ with respect to $\Phi_{t,\mu}$ is a free
module over $\C[[t]][\mu]$. For all $\mu\in\c^*$, the algebra
$\A_{t,\mu}$ specified to the fiber $O_\mu^k$ gives a
$G$-equivariant quantization of the KKS bracket on it.
\end{thm}
In particular, this theorem gives $G$-equivariant quantizations on
all semisimple orbits.

\section{$U_h(\g)$-equivariant quantization on orbits}
\label{sQGC}
\subsection{Extended reflection equation algebra and $(t,h)$ quantization
 on $gl(n)^*$}
We specialize our further considerations to the case
$\g=gl(n,\C)$. The $G$-module $\g^*$ is identified with
$\End(\C^n)$ using the trace pairing. The multiplication on matrix
units $\{e^i_j\}_{i,j=1}^n\subset \End(\C^n)$ is given by
$e^i_je^l_k=\delta^l_je^i_k$, where $\delta^l_j$ is the Kronecker
symbol. Let $S\in \End^{\tp 2}(\C^n)$ be the Hecke symmetry
related to the representation of $\U_h(\g)$ on $\C^n$. It
satisfies the relations \be{} \label{YBE}
S_{12}S_{23}S_{12}&=&S_{23}S_{12}S_{23},\\
S^2-(q-q^{-1})S&=&1\tp 1, \quad q=e^h. \label{Hecke} \ee{} The
symmetry $S$ is expressed through the image $R$ of the universal
R-matrix of $\U_h(\g)$ and the flip operator $P$ on $\C^n\tp
\C^n$, $S=PR$. Let $\{L^i_j\}\subset \End^*(\C^n)$ be the dual
basis to $\{e^j_i\}$. Consider an associative algebra, $\L_{t,h}$,
over $\C[[t,h]]$ generated by $L^i_j$ subject to the
quadratic-linear relations \be{} \label{ere_rel} S L_2S L_2- L_2S
L_2S=  qt \bigl(L_2S-SL_2\bigr). \ee{} Here,  $L$ is the matrix
$\sum_{i,j}L^i_je^j_i$ and $L_2=1\tp  L$. Note that at $h=0$ these
relations give the algebra $\L_t$ defined in Subsection 2.1.

Let $\Ru$ be the universal R-matrix of the quantum group
$\U_h(\g)$. Consider the element $\bar\QQ=qt/(q-q^{-1})(\QQ-1\ot
1)$, where $\QQ=\Ru_{21}\Ru$. It is easy to check that the map
$\L_{t,h}\to \U_h(\g)[[t]]$, $L^i_j\mapsto
L^i_j(\bar\QQ_1)\bar\QQ_2$ respects relations (\ref{ere_rel}). So
we have the homomorphism \be{} \label{multq}
\phi_{t,h}\colon\L_{t,h}\to \U_h(\g)[[t]]. \ee{}
\begin{propn}\label{doubleq}
The map $\phi_{t,h}$ is an embedding of free $\C[[t,h]]$-modules.
\end{propn}

\begin{proof} At the point $h=0$, the quotient of the tensor algebra
$\T\bigl(\End^*(\C^n)\bigr)$ by relations (\ref{ere_rel})
coincides with $\L_t$ and the map $\phi_{t,0}$ coincides with
(\ref{mult}), so the proposition is a corollary of  the following
lemma.
\end{proof}

\begin{lemma}\label{lem1}
Let $\Ec$ and $\F$ be vector spaces over $\C$ and $\la$ a set of
formal parameters. Let $\psi_\la\colon \Ec[[\la]]\to \F[[\la]]$,
$\psi_\la = \psi_0 \!\!\mod \!\la $, be a map of free
$\C[[\la]]$-modules and $\KK_\la$ a submodule in $\Ec[[\la]]$
satisfying the conditions: a) $\KK_\la\subset
\mathrm{ker}\:\psi_\la$, b) $\KK_0=\mathrm{ker}\:\psi_0$, where
$\KK_0=\{a\in \Ec|\exists f=a+\la b\in \KK_\la\}$. Then, the
$\C[[\la]]$-module $\Ec[[\la]]/\KK_\la$ is free and the map
$\Ec[[\la]]/\KK_\la\to \F[[\la]]$ is an embedding.
\end{lemma}
\begin{proof} Clear.
\end{proof}
\begin{corollary}
The algebra $\L_{t,h}$ is a $\U_h(\g)$-equivariant quantization on
$gl^*(n,\C)$.
\end{corollary}

\subsection{Double quantization on orbits in $gl^*(n,\C)$}
Let $\g=\End(\C^n)$. As a Cartan subalgebra in $\g$ we take
diagonal matrices. As a Levi subalgebra $\l$ we take the subspace
of matrices commuting with a diagonal matrix with eigenvalues
$\mu_1,...,\mu_{k}$. Then $\c$ consists of matrices commuting with
$\l$. Abusing notations, we will consider $(\mu_1,...,\mu_{k})$ as
coordinates in $\c^*\simeq\C^k$.

Since $\U_h(\g)\simeq \U(\g)[[h]]$ as associative algebras, there
is an extension,
$$\ff_{h,\mu}\colon\U_h(\g)\to \End(\V)[[h]],$$
of map (\ref{ff}). This algebra homomorphism is automatically
equivariant with respect to the adjoint actions of the Hopf
algebra $\U_h(\g)$. It is easy to check that the map
$\ff_{h,\mu/t}\circ\phi_{t,h}\colon \LL_{t,h}\to \End(\V)[[t,h]]$,
where $\phi_{t,h}$ is given by (\ref{multq}), defines a map of
polynomial sections in $\mu$, \be{}\label{Polsecq} \tilde
\Phi_{t,h,\mu}\colon\LL_{t,h}[\mu]\to \End(\V)[[t,h]][\mu]. \ee{}

\begin{thm}\label{thm2}
The image $\tilde \A_{t,h,\mu}$ of $\L_{t,h}[\mu]$ in
$\End(\V)[[t,h]][\mu]$ with respect to $\tilde \Phi_{t,h,\mu}$ is
flat over $\C[[t,h]][\mu]$. For a fixed $\mu\in\c^*$, $\tilde
\A_{t,h,\mu}$ gives a $(t,h)$-quantization on $\O^k_\mu$.
\end{thm}
\begin{proof}
The statement follows from Theorem \ref{thm1} and Lemma
\ref{lem1}.
\end{proof}

In the classical case, elements of $\O^k_\mu$ satisfy the matrix
equation \be{} \label{me}
(X-\mu_1)...(X-\mu_k)=X^{k}+\xi_{k-1}(\mu)X^{k-1}+...+\xi_0(\mu)=0,
\ee{} where coefficients $\xi_i$ are symmetric functions of $\mu$.
Entries of the matrix $X^m$, $m\in \N$, generate a
$\U(\g)$-submodule in $\tilde \A_{0,0,\mu}$ of
$\End^*(\C^n)$-type, which contains a one-dimensional trivial
module and a submodule of $sl(n,\C)$-type. Moreover, invariants
form a one dimensional submodule in $\tilde \A_{0,0,\mu}$, while
the $sl(n,\C)$-isotypic component has multiplicity $k-1$ and is
generated by $\bar X^{1},\ldots,\bar X^{k-1}$, where $\bar X^m$ is
the $sl(n,\C)$-type component in $X^m$.

Since $\tilde \A_{t,h,\mu}$ is a $\U_h(\g)$-equivariant
$\C[[t,h]][\mu]$-flat quotient of $\LL_{t,h}[\mu]$, we conclude
that $\tilde \A_{t,h,\mu}$ has the same decomposition into
irreducible components as $\tilde \A_{0,0,\mu}$, the classical
algebra of polynomial functions on $\O^k_\mu$, and (\ref{me})
extends to a matrix polynomial equation,
 \be{} \label{meq}
L^{k}+\xi_{k-1}(t,h,\mu)L^{k-1}+...+\xi_{0}(t,h,\mu)=0.
\ee{}
Let us fix $\mu$ such that $O^k_\mu$ is a semisimple orbit. Let
$U\subset \C^k$ be a neighborhood of $\mu$ such that $\O^k_{\mu'}$
is a semisimple orbit for $\mu'\in U$. Then, we can think of
elementary symmetric functions $\xi_i(\mu')$ of $\mu'$ as local
coordinates on $U$. Let $\A(U)$ denote the algebra of analytic
functions on $U$. Coefficients of the polynomial in the left hand
side of (\ref{meq}) define an algebra map $\psi\colon
\A(U)[[t,h]]\to \A(U)[[t,h]]$,
$(\xi_{k-1}(\mu'),...,\xi_0(\mu'),t,h)\mapsto
\bigl(\xi_{k-1}(t,h,\mu'),...,\xi_{0}(t,h,\mu'),t,h\bigr)$. Since
$\psi$ is identical at $(t,h)=(0,0)$, it is an algebra
automorphism. Therefore, $\psi$ defines a transformation of
parameters $(t,h,\mu')$ in the neighborhood $U$ of $\mu$. These
arguments lead to the following
\begin{propn}
\label{poleqorb} For every semisimple orbit ${\O^k_\mu}$, there
exists a $(t,h)$-quantization, $\A_{t,h,\mu}$, on ${\O^k_\mu}$ and
an epimorphism \be{} \label{map2} \Phi_{t,h,\mu}\colon\L_{t,h}\to
\A_{t,h,\mu} \ee{} factored through the ideal
$\J_{t,h,\mu}\subset\L_{t,h}$ generated by the relations \be{}
\label{meqq}
 (L-\mu_1)...(L-\mu_k)=L^k+\xi_{k-1}(\mu)L^{k-1}+...+\xi_0(\mu)=0.
\ee{}
\end{propn}
\begin{proof}
Let $\tilde\I_{t,h,\mu'}$ be the kernel of the map
$\tilde\Phi_{t,h,\mu'}$ from (\ref{Polsecq}), where $\mu'$ are
considered as local coordinates on the neighborhood $U$ of the
point $\mu$. Let $\I_{t,h,\mu}=\tilde \I_{\psi^{-1}(t,h,\mu)}$ and
$\A_{t,h,\mu}=\L_{t,h}/\I_{t,h,\mu}$. Since $\psi^{-1}$ transforms
the coefficients of polynomial (\ref{meq}) to those of polynomial
(\ref{meqq}), the natural map $\L_{t,h}\to\A_{t,h,\mu}$ satisfies
the proposition.
\end{proof}

\begin{remark}
\label{notminimal} Polynomial (\ref{meqq}) is minimal for the
semisimple orbit ${\O^k_\mu}$. Proposition \ref{poleqorb}
obviously holds for any polynomial divided by (\ref{meqq}).
\end{remark}
In the next subsection, we show that the quotient of the algebra
$\L_{t,h}$ by the ideal generated by entries of {\em any} matrix
polynomial is a flat deformation.

\subsection{Flatness of matrix polynomial relations}

Let $a=(a_{m-1},...,a_0)$ be coordinates in $\C^m$. Let us
consider the subvariety $M^m$ in $\g^*\times \C^m$ generated by
the matrix equation \be{}\label{meg} X^m+a_{m-1}X^{m-1}+...+a_0=0.
\ee{} It is clear that if the polynomial
$f_a(x)=x^m+a_{m-1}x^{m-1}+...+a_0$ has no multiple roots, the
fiber $M^m_a$ is the union of semisimple orbits passing through
diagonal matrices with eigenvalues from the set $\{\la_i\}$,
$i=1,...,m$, of simple roots of $f(x)$. It is easy to see that the
variety $M^m$ is flat over $\C^m$.

\begin{thm}\label{thm3}
Let $\J_a$ be the ideal in the ERE algebra $\LL_{t,h}$ generated
by entries of a matrix polynomial relation,
\be{}\label{megg}
f_a(L)=L^m+a_{m-1}L^{m-1}+...+a_0=0. \ee{} Then the quotient
algebra $\LL_{t,h}/\J_a$ is flat over $\C[[t,h]][a]$ and for any
$a\in\C^m$ gives a $(t,h)$-quantization on $M^m_a$.
\end{thm}
\begin{proof}
First, we suppose that $f_a(x)$ has only simple roots
$(\la_1,\ldots,\la_m)$. Let the fiber $M^k_a$ be the union of
semisimple orbits $O_\ell$, $\ell=1,...,K$.
As pointed out in Remark \ref{notminimal}, for every $\ell$ there
exists an epimorphism $\Phi^\ell_{t,h,}\colon \L_{t,h}/\J_a\to
\A_{t,h}^\ell$, where $\A_{t,h}^\ell$ is a $(t,h)$-quantization on
$\O_\ell$. Consider the algebra homomorphism \be{}\label{ahom}
\oplus_{\ell=1}^K\Phi^\ell_{t,h}\colon\L_{t,h}/\J_a\to\oplus_{\ell=1}^K\A_{t
,h}^\ell. \ee{} This map is an isomorphism at the classical point
$(t,h)=(0,0)$. Therefore, as follows from Lemma \ref{lem1}, it is
an isomorphism in the quantum case, too. This proves flatness of
$\LL_{t,h}/\J_a$ over neighborhoods of all $a\in\C^m$
corresponding to polynomials with simple roots.

Let $f_a(x)$ be an arbitrary polynomial. Then, it is a limit of
polynomials with simple roots. Now, the flatness of
$\LL_{t,h}/\J_a$ over $\C[[t,h]][a]$  follows from the above part
of the proof and from the fact that the quotient algebra
$\L_{0,0}/\J_a$ is flat over $\C[a]$.
\end{proof}

The orbits in $M^m_a$ correspond to points of the spectrum of the
subalgebra of invariants in the function algebra on $M^m_a$. In
order to pass to the quantization on a particular orbit, we should
study the subalgebra of invariants in  $\L_{t,h}$ restricted to
$M^m_a$. That is done in the next section.

\section{Algebra of invariants on $M^k$}
\label{sInv}
\subsection{Center of the extended reflection equation algebra}
The subalgebra $\S_{t,h}\subset \L_{t,h}$ of
$\U_h\bigl(gl(n,\C)\bigr)$-invariants coincides with the center of
$\L_{t,h}$. It is generated by elements $s_1,\ldots,s_n$, where
$s_m=\Tr_q( L^m)$, $m\in \N$. The quantum trace $\Tr_q$ of a
matrix $A$ is defined as, \cite{FRT},
\be{}
\Tr_q(A)=\Tr(DA),\quad q=e^{h},
\ee{}
with the weight matrix $D$ normalized by the condition
$(\id\tp\Tr_q)(S)=q \:1$. In this normalization, $\Tr_q(1)=\hat
n$. We use the notation $\hat m$ for quantum integers, $\hat
m=\frac{1-q^{-2m}}{1-q^{-2}}$, $m\in \Z$.

The algebra $\S_{t,h}$ is isomorphic to
$\C[[t,h]][\si_1,\ldots,\si_n]$, where the elements $\si_i$ are
coefficients of the characteristic polynomial in the quantum
Cayley-Hamilton identity \be{} \label{Cayley-Hamilton}
 L^n - \si_1  L^{n-1}+\ldots + (-1)^n\si_n=0
\ee{} which follows from the same representation theory arguments
as in the proof of (\ref{meq}).

We substantially rely on relations between generators $\{\si_i\}$
and $\{s_i\}$ calculated by Pyatov and Saponov, \cite{PS}, for the
algebra $\S_{h}=\S_{0,h}$. In our normalization of the quantum
trace, their result is formulated as follows.
\begin{thm}[\cite{PS}]
\label{PSthm} Elements $\{s_i\}_{i=1}^\infty$ and
$\{\si_i\}_{i=1}^\infty$, where $\si_i=0$ for $i>n$, of the
algebra $\S_h$ satisfy the quantum Newton identities \be{}
\label{Newton} \hat m  \si_m - s_1\si_{m-1}+\ldots+(-1)^ms_{m}=0,
\quad m=1,2,\ldots. \ee{}
\end{thm}
Given a sequence $\vec x=\{x_m\}_{m=0}^\infty$ with values in a
$\C[\mu]$-module, where $\mu=(\mu_1,...,\mu_k)\in\C^k$, $k\leq n$,
let us define the sequence $\{r_m(\vec x)\}_{m=k}^\infty$, \be{}
\label{ri} r_m(\vec x)=x_{m}- \si_1(\mu) x_{m-1}+\ldots+(-1)^k
\si_k(\mu) x_{m-k}, \ee{} where $\si_i(\mu)$ are the elementary
symmetric functions in $\mu$, $\si_i(\mu)=\sum_{1\leq
j_1<\ldots<j_i\leq k} \mu_{j_1}\ldots\mu_{j_i}$. The following
proposition describes the intersection of $\S_{t,h}$ with the
ideal $\J_{t,h,\mu}\subset\L_{t,h}$ generated by relations
(\ref{meqq}).
\begin{propn}
\label{recurrent} Let $\vec s=\{s_m\}_{m=0}^\infty$ be the
sequence of q-traces, $s_m=\Tr_q(L^m)$. Then, the set $\{r_m(\vec
s)\}_{m=k}^\infty$ generates the ideal $\J^\S_{t,h,\mu}=\S_{t,h}
\cap \J_{t,h,\mu}$ in $\S_{t,h}$.
\end{propn}
\begin{proof}
Denote by $p$ the polynomial on the left-hand side of (\ref{meqq})
whose entries $p^i_j$ generate the ideal $\J_{t,h,\mu}$. Clearly
$r_m(\vec s)=\Tr_q\bigl(L^{m-k}p\bigr)$, so the elements $r_m(\vec
s)$ belong to $\J^\S_{t,h,\mu}$ for all $m\geq k$. Let us prove
that they generate $\J_{t,h,\mu}^\S$. Suppose $y\in \J_{t,h,\mu}$;
then it is representable as $y=\sum_{i,j=1}^n y^i_jp^j_i$, where
$\{y^i_j\}\subset \L_{t,h}$. If $y\in \S_{t,h}$, the elements
$y^i_j$ generate a module of $\End(\C^n)$-type. But the
$\End(\C^n)$-type component in $\L_{t,h}$ is an $\S_{t,h}$-module
spanned by entries of the matrix powers $L^l$, $l=0,\ldots,n-1$.
Therefore $y^i_j$ are entries of a polynomial $\sum_{l=0}^{n-1}
y_l\: L^l$ with coefficients $y_l\in \S_{t,h}$, so
$y=\sum_{l=0}^{n-1} y_l\:r_{l+k}(\vec s)$.
\end{proof}
\subsection{Special polynomials}
\label{spec_pol}
Characters of the algebra $\S_{t,h}$ restricted
to $M^k$ are described by means of special polynomials, which are introduced in
this subsection.

Let us define the set $\{n\!:\!k\}\subset \Z^k$ of $k$-tuples
${\bf n}=(n_1,\ldots,n_k)$ such that $0\leq n_i$ and
$n_1+\ldots+n_k=n$; the subset $\{n\!:\!k\}_+\subset \{n\!:\!k\}$
consists of ${\bf n}$ with all $n_i>0$. We denote by $|n\!:\!k|$
the number of elements in $\{n\!:\!k\}$.

Given ${\bf n}\in \{n\!:\!k\}$ let us introduce polynomials in
$\mu$ of degree $m=0,1,\ldots $ setting \be{} \label{tdef}
\vartheta_m({\bf n},q^{-2},\mu) &=& \sum_{\ell=1}^k
(1-q^{-2})^{\ell-1} \sum_{{\bf d}\in \{m:\ell\}_+}\; \sum_{1\leq
i_1<\ldots <i_\ell\leq k} \hat n_{i_1}\ldots \hat n_{i_\ell}
\:\mu^{d_1}_{i_1} \ldots\mu^{d_\ell}_{i_\ell} \ee{} for $m>0$  and
$\vartheta_0({\bf n},q^{-2},\mu)=\hat n$. At the classical point
$q=1$, the polynomial  $\vartheta_m({\bf n},q^{-2},\mu)$ is equal
to $\sum_{i=1}^kn_i\mu^m_i$, i.e., the trace $\Tr(A^m)$ of a
matrix $A$ with eigenvalues $(\mu_1,\ldots,\mu_k)$ of
multiplicities $(n_1,\ldots,n_k)$.

The following proposition will be important for our consideration.
\begin{propn}
The polynomials $\vartheta_m({\bf n},q^{-2},\mu)$, $m\geq 0$, can
be represented as \be{} \vartheta_m({\bf
n},q^{-2},\mu)&=&\sum_{j=1}^k C_j({\bf n},q^{-2},\mu)\mu_j^m,
\label{rational} \ee{} \be{} \mbox{where}\quad C_j({\bf
n},q^{-2},\mu)&=& \hat n_j+\hat n_j \sum_{\ell=1}^{k-1}
(1-q^{-2})^{\ell} \sum_{1\leq i_1<\ldots <i_\ell\leq k\atop
i_1,\ldots, i_\ell\not=j} \frac{\hat n_{i_1}
\mu_{i_1}}{\mu_{j}-\mu_{i_1}}\ldots \frac{\hat n_{i_\ell}
\mu_{i_\ell}}{\mu_{j}-\mu_{i_\ell}} \:. \ee{}
\end{propn}
\begin{proof}
It follows from definition (\ref{tdef}) that the polynomials
$\vartheta_m({\bf n},q^{-2},\mu)$ satisfy the identity \be{}
\label{n_n-1} \vartheta_m({\bf n},q^{-2},\mu) =\vartheta_m({\bf
n'},q^{-2},\mu')+(1-q^{-2})\hspace{2pt}\hat n_k
\sum_{i=1}^{m-1}\vartheta_{m-i}({\bf
n'},q^{-2},\mu')\mu_k^i + \hat n_k \mu_k^m, \ee{} where
$\mu'=(\mu_1,\ldots,\mu_{k-1})$ and ${\bf
n}'=(n_1,\ldots,n_{k-1})$. Equation (\ref{n_n-1}) allows to apply
induction on $k$.
\end{proof}

Using representation (\ref{rational}), let us introduce  the
functions
$$
\vartheta_m({\bf n},q^{-2},\mu,t)=\sum_{j=1}^k C_j\Bigl({\bf
n},q^{-2},\mu+\frac{t}{1-q^{-2}}\Bigr)\mu_j^m,
$$
where we put
$\mu+\frac{t}{1-q^{-2}}=(\mu_1+\frac{t}{1-q^{-2}},\ldots,\mu_k+\frac{t}{1-q^
{-2}})$. Although the coefficients $C_j({\bf n},q^{-2},\mu)$ are
rational, the functions $\vartheta_m({\bf n},q^{-2},\mu,t)$ are,
in fact, polynomials in $t$, $q^{-2}$, and $\mu$. Evidently,
$\vartheta_m({\bf n},q^{-2},\mu)=\vartheta_m({\bf
n},q^{-2},\mu,t)|_{t=0}$.

\subsection{Characters of the subalgebra of invariants in $\L_{0,h}$}
Let us reserve the notation $\la=\mu$ for the case $k=n$ and
consider polynomials $\vartheta_m(q^{-2},\la)$; by definition,
\be{} \label{tm} \vartheta_m(q^{-2},\la)=\vartheta_m({\bf
n},q^{-2},\la)|_{{\bf n}=(1,\ldots,1)} =\sum_{\ell=1}^{n}
(1-q^{-2})^{\ell-1} \;\sum_{{\bf d}\in \{m:\ell\}_+}\; \sum_{1\leq
j_1 < \ldots < j_\ell\leq n}
\la^{d_1}_{j_1}\ldots\la^{d_\ell}_{j_\ell}. \ee{}
\begin{lemma}
\label{K=1} Substitution $\la_i=q^{-2(i-1)}$, $i=1,\ldots, n$, to
$\vartheta_m(q^{-2},\la)$ returns $\hat n $.
\end{lemma}
\begin{proof}
We apply induction on $n$. For $n=1$ we obviously have
$\vartheta_m(q^{-2},\la)=1$. Suppose the theorem holds for $n=l$.
Using the induction assumption and identity (\ref{n_n-1}) for the
case ${\bf n}=(1,\ldots,1)$ we find, for $n=l+1$, \be{}
\vartheta_m(q^{-2},\la)|_{\la=(1,q^{-2},\ldots,q^{-2l})} &=& \hat
l  + ( q^{-2l})^m + (1-q^{-2}) \hat l\sum_{d=1}^{m-1} (q^{-2l})^d
\nn\\&=& \hat l + (q^{-2l})^m + (1-q^{-2})\hat l(q^{-2l})
\frac{1-(q^{-2l})^{m-1}}{1- q^{-2l}}\:. \nn \ee{} After elementary
transformations the last expression is brought to $\widehat {l+1}
$.
\end{proof}
Given ${\bf n}=(n_1,\ldots,n_k)\in \{n\!:\!k\}$ let us consider
non-intersecting intervals $\Ind_i\subset \{1,\ldots, n\}$,
$i=1,\ldots,k$, some of them may be empty, satisfying the
following requirements: 1) if $a\in \Ind_i$ and $b \in \Ind_j$, then
$a<b$  whenever $i<j$, 2) $\#\Ind_i=n_i$. Clearly, $\Ind_1\cup
\ldots \cup\Ind_k=\{1,\ldots, n\}$. Let us introduce double
indexing $\la_{j,i}$, $j=1,\ldots,k$, $i=1,\ldots,n_j$, of the
indeterminates $\la_1,\ldots,\la_n$: if $m\in \Ind_j$ and $i$ is
the relative position of $m$ within the interval $\Ind_j$, then we
identify $\la_{j,i} =\la_m$. Double indices are ordered by the
lexicographic ordering, which coincides with that induced from
$\N$. We shall consider the  vectors $\bar
\la_j=(\la_{j,1},\ldots,\la_{j,n_j})$ and polynomials
$\vartheta_m(q^{-2},\bar\la_j)$, $j=1,\ldots,k$.
\begin{propn}
\label{q-substitution} Substitution $\la_{j,i}=\mu_j q^{-2(i-1)}$
to  $\vartheta_m(q^{-2},\la)$, $m\in \N$, gives
$\vartheta_{m}({\bf n},q^{-2},\mu)$.
\end{propn}
\begin{proof}
The case $k=1$ follows from Lemma \ref{K=1}, because the
polynomials $\vartheta_m(q^{-2},\la)$, $m\in \N$, are homogeneous
in $\la$ of degree $m$. The general situation is reduced to the
case $k=1$ if one observes that, upon rearranging summation in
(\ref{tm}), the polynomial $\vartheta_m(q,\la)$ may be rewritten
as \be{} \label{ind} \vartheta_m(q^{-2},\la)= \sum_{\ell=1}^k
(1-q^{-2})^{\ell-1} \:\sum_{{\bf d}\in \{m:\ell\}_+}\: \sum_{1\leq
j_1 < \ldots < j_\ell\leq k} \vartheta_{d_1}(q^{-2},\bar
\la_{j_1})\ldots \vartheta_{d_\ell}(q^{-2},\bar \la_{j_\ell}),
\ee{} Now we apply Lemma \ref{K=1} to each factor
$\vartheta_{d_1}(q^{-2},\bar \la_{j_1}), \ldots
,\vartheta_{d_\ell}(q^{-2},\bar \la_{j_\ell})$ on
 the right-hand side of
(\ref{ind}).
\end{proof}

Recall that the subalgebra of invariants $\S_{h}=\S_{0,h}$ of the
quadratic RE algebra $\L_h=\L_{0,h}$ is a polynomial algebra in
$n$ variables $\si_i$, $i=1\ldots,n$, which are the coefficients
of the Cayley-Hamilton identity (\ref{Cayley-Hamilton}). Let us
assign to $\la\in \C^n[[h]]$ a character
$\chi^\la\colon\S_{h}\to \C[[h]]$
setting on the generators
$\chi^\la(\si_m)=\sum_{i_1<\ldots<i_m}\la_{i_1}\ldots
\la_{i_m}=\si_m(\la)$, the elementary symmetric polynomial of
degree $m\in\{1,\ldots,n\}$. The correspondence $\la\to \chi^\la$
defines an epimorphism from $\C^n[[h]]$ to $\Char\: \S_h$, the set
of characters of $\S_h$.
\begin{propn}
One has $\chi^\la(s_m)=\vartheta_m(q^{-2},\la)$, where
$s_m=\Tr_q(L^m)$, $m\in \N$. \label{q-traces}
\end{propn}
\begin{proof}
Assuming $\si_m=0$ and $\si_m(\la)=0$ for $m>n$ we must check that
the substitution $\si_m\to\si_m(\la)$, $s_m\to
\vartheta_m(q^{-2},\la)$ satisfies Newton identities
(\ref{Newton}). Since representation (\ref{rational}) is valid for
$\vartheta_m(q^{-2},\la)=\vartheta_m({\bf n},q^{-2},\la)|_{{\bf
n}=(1,\ldots,1)}$, it is enough to check identities (\ref{Newton})
for $m\leq n$ only. Let $i$ and $j$ be non-negative integer
numbers and integers $d_1,\ldots,d_j>1$ be such that $i+d_1+\ldots
+d_{j}=m\leq n$. Consider the coefficient before the monomial
$p=p'p''$ in equality (\ref{Newton}), where $p'=\la_1 \ldots
\la_i$ and $p''=\la^{d_1}_{i+1} \ldots \la^{d_j}_{i+j}$. Let us
show that this coefficient is equal to zero. This will prove the
statement because all the polynomials involved are symmetric in
$\la$. Let us compute contributions coming, say, from the term
$(-1)^l\si_l(\la) s_{m-l}(\la)$, where $l\geq0$ (we assume
$\si_0=1$). We should consider all the monomials $\frac{p}{r}$,
where $r=\la_{\al_1}\ldots \la_{\al_l}$, $\al_1<\ldots<\al_l\leq
i+j$, enters $\si_l(\la)$ and take the coefficients before
$\frac{p}{r}$ entering $s_{m-l}(\la)$, with the sign $(-1)^l$. We
can represent $r$ as the product $r=r'r''$, where $r'$ contains
only $\la_l$, $1\leq l\leq i$, while $r''$ involves only $\la_l$,
$i<l\leq i+j$; so $r'$ and $r''$ divide $p'$ and $p''$,
respectively. Observe that the coefficients before $\frac{p}{r}$
in $s_m(\la)$ do not depend on $r''$ and are equal to
$(1-q^{-2})^{i+j-\deg(r')-1}$.

Consider separately two cases: $p''=1$ and $p''\not=1$. If
$p''=1$, then the coefficient before the  term $\la_1\ldots \la_m$
entering (\ref{Newton}) is \be{}
 \hat m +  \sum_{l=1}^{m} (-1)^{l} {l \choose m}\:(1-q^{-2})^{l-1} = \hat
m + \frac{\bigl(1-(1-q^{-2})\bigr)^m-1}{(1-q^{-2})}=0. \ee{}

If $p''\not =1$, the first term in  (\ref{Newton}) gives no
contribution to the coefficient, which we find to be
$$
\sum_{r'}(1-q^{-2})^{i+j-\deg(r')-1}\sum_{r''} (-1)^{\deg(r)}=
\sum_{r'}(1-q^{-2})^{i+j-1}\sum_{l=0}^j (-1)^l {l\choose j}=0.
$$
\end{proof}

\subsection{Algebra of invariants on $M^k$ and its restriction to orbits}
In this section, we study the image of the subalgebra of
invariants $\S_{t,h}\subset \L_{t,h}$ in the quotient of
$\L_{t,h}$ by the ideal $\J_{t,h,\mu}$ generated by relations
(\ref{meqq}). To analyze its structure and compute its characters,
we use polynomials $\vartheta_m({\bf n},q^{-2},\mu)$ and
$\vartheta_m({\bf n},q^{-2},\mu,t)$ introduced in Subsection \ref{spec_pol}.
\begin{thm}
\label{quantum} The quotient  of the algebra $\S_{t,h}$ by the
ideal $\J^\S_{t,h,\mu}=\S_{t,h}\cap\J_{t,h,\mu}$ is generated by
the images of the elements $\{s_1,\ldots,s_{k-1}\}$ in
$\S_{t,h}/\J^\S_{t,h,\mu}$. For any ${\bf n}\in\{n\!:\!k\}$, the
formula \be{} \label{char_q}
 \chi_{t,h}({\bf n})(s_m)= \vartheta_m({\bf n},q^{-2},\mu,t),
 \quad m\in \N,
\ee{} defines a character $\chi_{t,h}({\bf n})$ of the algebra
$\S_{t,h}/\J^\S_{t,h,\mu}$. The correspondence ${\bf n}\to
\chi_{t,h}({\bf n})$ is a bijection between $\{n\!:\!k\}$ and the
set of characters $\Char\:\S_{t,h}/\J^\S_{t,h,\mu}$. The algebra
$S_{t,h}/\J^\S_{t,h,\mu}$ is a free $\C[[t,h]]$-module  of rank
$|n\!:\!k|$. At generic $\mu$, it is a direct sum of  $|n\!:\!k|$
copies of  $\C[[t,h]]$.
\end{thm}
\begin{proof}
First of all, the theorem holds true in the classical case
$(t,h)=(0,0)$. Clearly the images of the  elements
$\{s_1,\ldots,s_{k-1}\}$ generate $\S_{t,h}/\J^\S_{t,h,\mu}$
because other traces are related to them via recurrent relations,
according to Proposition \ref{recurrent}.

Let us describe characters of the algebra
$\S_{t,h}/\J^\S_{t,h,\mu}$. Suppose first $t=0$. An element $\chi$
from $\Char\:\S_{0,h}/\J^\S_{0,h,\mu}$ is also a character for
$\S_{0,h}$, so we can assume $\chi=\chi^\la$, for some $\la\in
\C^n[[h]]$ depending on $\mu$. We have
$\chi(s_m)=\vartheta_m(q^{-2},\la)$ according to Proposition
\ref{q-traces}. For $\chi^\la$ to define a character of the
quotient algebra $\S_{0,h}/\J^\S_{0,h,\mu}$, the sequence
$\chi(\vec s)\subset \C[[h]]$, must satisfy the recurrent
identities $r_m\bigl(\chi(\vec s)\bigr)=0$, $m\geq k$, by
Proposition \ref{recurrent}. By Proposition \ref{q-substitution},
for any ${\bf n}\in \{n\!:\!k\}$ the polynomials $\vartheta_m({\bf
n},q^{-2},\mu)$ are obtained by a substitution to
$\vartheta_m(q^{-2},\la)$. The sequence $\vec a=\{\vartheta_m({\bf
n},q^{-2},\mu)\}_{m=0}^\infty$ satisfies the recurrent relations
$r_m(\vec a)=0$, $m\geq k$, as follows from representation
(\ref{rational}).
 Therefore, for any ${\bf n}\in \{n\!:\!k\}$ there exists  an element
$\chi_h({\bf n})\in\Char\:\S_{0,h}/\J^\S_{0,h,\mu}$ such that
$\chi_h({\bf n})(s_m)=\vartheta_m({\bf n},q^{-2},\mu)$, for $m\geq
0$.

  As follows from their definitions, the ERE and RE algebras are defined
over $\C[q,q^{-1}]$, $q=e^h$, and all the $\C[[h]]$-algebras under
consideration are, in fact, completions of the corresponding
$\C[q,q^{-1}]$-algebras at the point $q=1$. Since the functions
$\vartheta_m({\bf n},q^{-2},\mu)$ are polynomials in $q^{-2}$, the
character $\chi_h({\bf n})$ is the extension of a character
$\chi_q({\bf n}) \in \Char \S_{0,q}/\J^\S_{0,q,\mu}$.

Let us prove formula (\ref{char_q}) for the case $t\not=0$. Hecke
condition (\ref{Hecke}) and relations (\ref{ere_rel}) imply that
the shift $L\mapsto L-\frac{t}{1-q^{-2}}$, $\mu \mapsto
\mu-\frac{t}{1-q^{-2}}$ induces a homomorphism
$\varrho_{t,q}\colon \L_{t,q}\to \L_{0,q}[t]$ (here we assume that
the ring $\C[q,q^{-1}]$ is extended to the field $\C(q)$). Clearly
$\varrho_{t,q}(\J_{t,q,\mu})$ coincides with  $\J_{0,q,\mu}$
because equation (\ref{meqq}) is invariant under this
transformation. Therefore, the map $\varrho_{t,q}$ defines a
homomorphism $\bar \varrho_{t,q}\colon\S_{t,q}/\J^\S_{t,q,\mu} \to
\S_{0,q}/\J^\S_{0,q,\mu}$. Denote by  $\chi_{t,q}({\bf n})$ the
composition $\chi_q({\bf n})\circ\bar \varrho_{t,q}$, which is a
character of $\S_{t,q}/\J^\S_{t,q,\mu}$. It is easy to prove by
induction on $m$, using representation (\ref{rational}), that
$\chi_{q,t}({\bf n})$ evaluated on $s_m$ gives $\vartheta_m({\bf
n},q^{-2},\mu,t)$. Since the functions  $\vartheta_m({\bf
n},q^{-2},\mu,t)$ are polynomials in $q^{-2}$, we can extend
$\chi_{t,q}({\bf n})$  to a character $\chi_{t,h}({\bf n})$ of
$\S_{t,h}/\J^\S_{t,h,\mu}$ over $\C[[h]]$.

Thus we have constructed a deformation $\chi_{t,h}({\bf n})$ for each
classical character $\chi({\bf n})$. We have the inequality
$$|n\!:\!k|=\#\Char\:\S/\J^\S_\mu\leq\#\Char\:\S_{t,h}/\J^\S_{t,h,\mu}\leq
\mathrm{rk}\:\S_{t,h}/\J^\S_{t,h,\mu}\leq N
\leq\mathrm{rk}\:\S/\J^\S_\mu =|n\!:\!k|,$$ where $\S=\S_{0,0}$
and $\J_{\mu}=\J^\S_{0,0,\mu}$ are the classical algebras and $N$
the number of generators of the $\C[[t,h]]$-module
$\S_{t,h}/\J^\S_{t,h,\mu}$. This inequality shows that
$\S_{t,h}/\J^\S_{t,h,\mu}$ is a free module over $\C[[t,h]]$ of
rank $|n\!:\!k|$.

Let us fix $\mu$ such that $\mu_i\not=\mu_j$ for  $i\not=j$. Since
at the classical point $(t,h)=(0,0)$ the algebra map
$$
\oplus_{{\bf n}\in \{n:k\}}\chi_{t,h}({\bf n})\colon
\S_{t,h}/\J^\S_{t,h,\mu}\to \oplus_{{\bf n}\in \{n:k\}}\C[[t,h]]
$$
is an isomorphism, it is an isomorphism of $\C[[t,h]]$-algebras.
\end{proof}
\section{Quantum  orbits}
\label{sQO} In this section, we present explicit formulas of  the
$\U_h\bigl(gl(n,\C)\bigr)$-equivariant quantization on orbits
using characters of the algebra $\S_{t,h}/\J_{t,h,\mu}^\S$
calculated in the previous section.
\begin{thm}
\label{final} For any element ${\bf n}=(n_1,\ldots,n_k)\in
\{n\!:\!k\}$, the quotient of the ERE algebra (\ref{ere_rel}) by
the relations \be{} \label{PR} (L-\mu_1)\ldots (L-\mu_k)&=&0,
\\
\label{TrR} \Tr_q(L^m)&=&\vartheta_m({\bf n},q^{-2},\mu,t), \quad
m=1,\ldots,k-1, \ee{} is a flat $\C[[t,h]]$-algebra when
$\mu_i\not=\mu_j$ for $i\not=j$. It is a quantization on the orbit
of semisimple matrices with eigenvalues $(\mu_1,\ldots\mu_k)$ of
multiplicities $(n_1,\ldots,n_k)$.
\end{thm}
\begin{proof}
By Theorem \ref{thm3}, the quotient algebra
$\L_{t,h}/\J_{t,h,\mu}$ defined by relations (\ref{PR}) is a
quantization on $M^k$. The algebra of invariants on a semisimple
orbit is one-dimensional thus corresponding to a character of
$\S_{t,h}/\J^\S_{t,h,\mu}$. At the classical point, characters are
$\chi({\bf n})(s_m)=\sum_{i=1}^k n_i\mu^m_i$, for ${\bf n}\in
\{n\!:\!k\}$. Their deformations are given by formula
(\ref{char_q}). So relations (\ref{TrR}) should hold on the orbit
of semisimple matrices with eigenvalues $(\mu_1,\ldots\mu_k)$ of
multiplicities $(n_1,\ldots,n_k)$. Because relations  (\ref{PR})
and  (\ref{TrR}) define this orbit at $(t,h)=(0,0)$, they define a
$(t,h)$-quantization on it.
\end{proof}
Since functions $\vartheta_m({\bf n},q^{-2},\mu,t)$ are
polynomials in $q^{-2}$, $t$, and $\mu$, we obtain the following
\begin{corollary}
At $q=1$, formulas (\ref{PR}) and (\ref{TrR}) give $G$-equivariant
quantizations of the Kostant-Kirillov-Souriau brackets on
semisimple orbits.
\end{corollary}

\begin{remark}
One can consider the ERE algebra associated with arbitrary Hecke
symmetry $S$ and define "quantum orbits" by formulas (\ref{PR})
and (\ref{TrR}). Equations (\ref{TrR}) are forced by condition
(\ref{PR}) and have the same form, since Newton identities
(\ref{Newton}) hold for any Hecke symmetry,  \cite{IOP}. The
question is about the module structure (supply of "functions") of
the orbits defined in this way. In our next paper we will prove that the supply of functions on
an "orbit" for any even Hecke symmetry is rich enough.

\end{remark}

\end{document}